\documentclass[11pt,a4paper]{article}
\usepackage{amsmath, amsfonts, amsthm, amssymb, graphicx}

\begin{document}

\title
{Counterexamples to a conjecture of Merker on 3-connected cubic planar graphs with a large cycle spectrum gap}
\author{
{\sc Carol T. ZAMFIRESCU\footnote{Department of Applied Mathematics, Computer Science and Statistics, Ghent University, Krijgslaan 281 - S9, 9000 Ghent, Belgium and Department of Mathematics, Babe\c{s}-Bolyai University, Cluj-Napoca, Roumania; e-mail address: \emph{czamfirescu@gmail.com}}}}
\date{}

\maketitle
\begin{center}
\vspace{12mm}
\begin{minipage}{125mm}
{\bf Abstract.} Merker conjectured that if $k \ge 2$ is an integer and $G$ a 3-connected cubic planar graph of circumference at least $k$, then the set of cycle lengths of $G$ must contain at least one element of the interval $[k, 2k+2]$. We here prove that for every even integer $k \ge 6$ there is an infinite family of counterexamples.

\bigskip

{\bf Key words.} Cycles; Cycle spectrum; 3-connected; Cubic; Planar graphs

\bigskip

\textbf{MSC 2020.} 05C38, 05C10

\end{minipage}
\end{center}

\vspace{12mm}

\section{Introduction}

For a graph $G$, we denote by ${\cal C}(G)$ the set of lengths of cycles in $G$, i.e.\ its \emph{cycle spectrum}. The \emph{circumference} of $G$ is the length of a longest cycle in $G$. Merker~\cite{Me21} recently proved that for any non-negative integer $k$ every 3-connected cubic planar graph $G$ of circumference at least $k$ satisfies ${\cal C}(G) \cap [k, 2k+9] \ne \emptyset$. He conjectured that for any integer $k \ge 2$ and any 3-connected cubic planar graph $G$ of circumference at least $k$, we have ${\cal C}(G) \cap [k, 2k+2] \ne \emptyset$. We shall abbreviate this conjecture of Merker with ($\dagger$).

By Euler's formula, every cubic plane graph contains a face of length 3, 4, or 5, so ($\dagger$) holds for $k \in \{ 2, 3 \}$. Suppose ($\dagger$) is untrue for $k = 5$. Then there exists a 3-connected cubic plane graph $G$ of circumference at least~$5$ with ${\cal C}(G) \cap [5, 12] = \emptyset$. Any 3- or 4-cycle in $G$ must be the boundary of a face of $G$, and any two faces in $G$ of size 3 or 4 are disjoint since $5 \notin {\cal C}(G)$ and $6 \notin {\cal C}(G)$. We contract every triangle and every quadrilateral of $G$ to a vertex and obtain the graph $G'$. If we exclude 3- and 4-cycles, cycles in $G$ have length at least 13, so $G'$ is a planar 3-connected graph with no cycle of length less than 7 (as on any $\ell$-cycle $C$ of $G$, $C$ shares at most $\lfloor \ell/2 \rfloor$ edges with a 3- or 4-cycle), a contradiction. The argument for $k = 4$ is very similar (and simpler). This yields that ($\dagger$) holds for $k \in \{ 2, 3, 4, 5 \}$. However, we now show that for any even integer $k \ge 6$ there is an infinite family of counterexamples to ($\dagger$).


\section{Result}

\noindent \textbf{Theorem.} \emph{For any even integer $k \ge 6$ there exists an infinite family of $3$-connected cubic planar graphs of circumference at least $k$ whose cycle spectrum contains no element of $[k,2k+2]$.}

\bigskip

\noindent \emph{Proof.} Consider the graph $H$ depicted in Fig.~1. Its left-most and right-most parts should be identified in the obvious way, where the boundary cycles of the two faces incident only with pentagons (top and bottom of Fig.~1) may have any length of at least $2r + 8$ (this yields the advertised infinite family). The vertices of $H$ are either black or white, as illustrated in Fig.~1.

\begin{center}
\includegraphics[height=70mm]{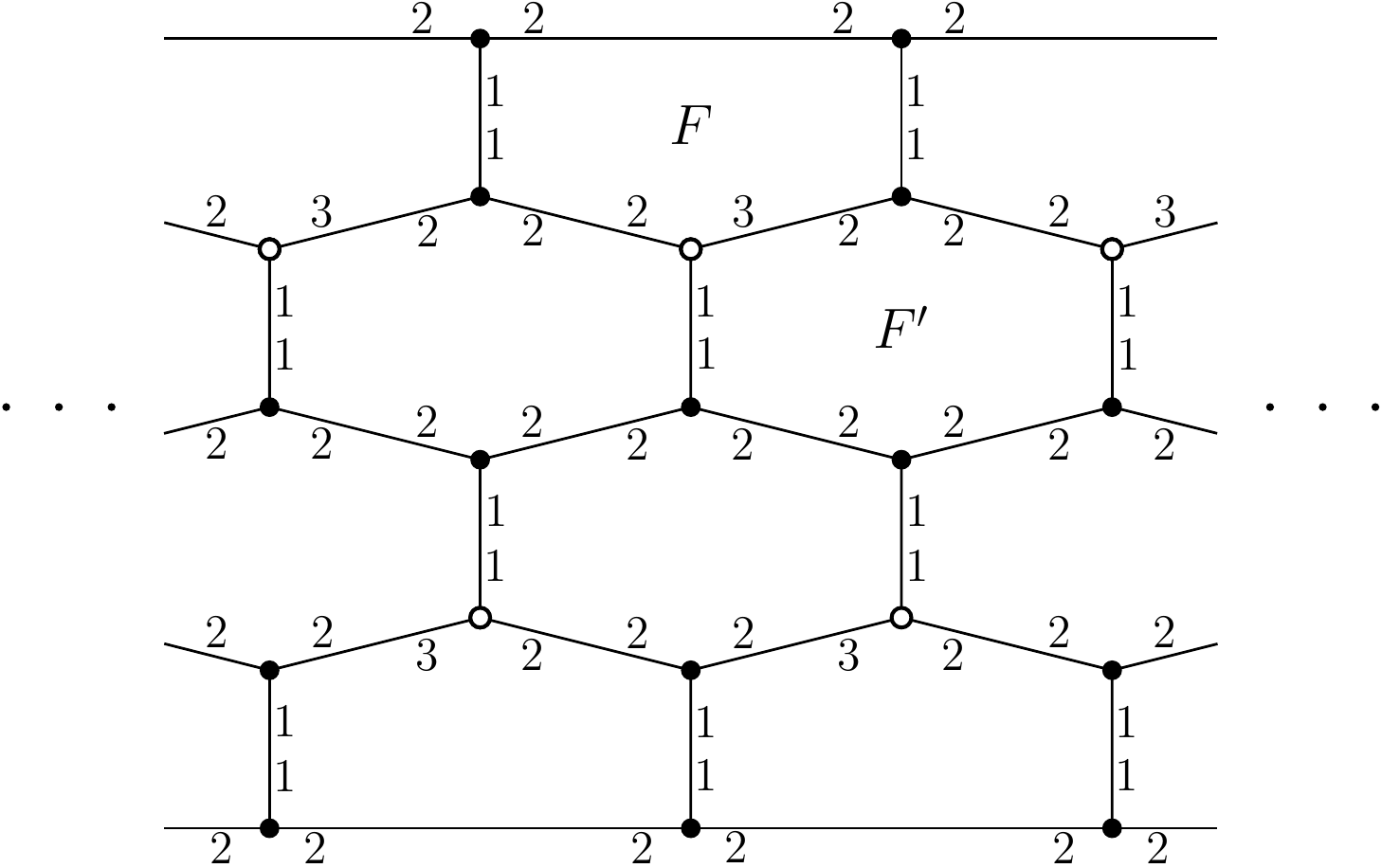}\\[1mm]
Figure 1: The graph $H$.
\end{center}

Consider the operations $A$ and $B$ defined in Fig.~2. We shall call a \emph{rung} any edge depicted as a horizontal line-segment in Fig.~2. In each operation, we replace a cubic vertex with the plane graph $A_r$ and $B_r$ (in which we ignore the three dangling edges), respectively, where $r$ denotes the number of rungs.

\begin{center}
\includegraphics[height=36mm]{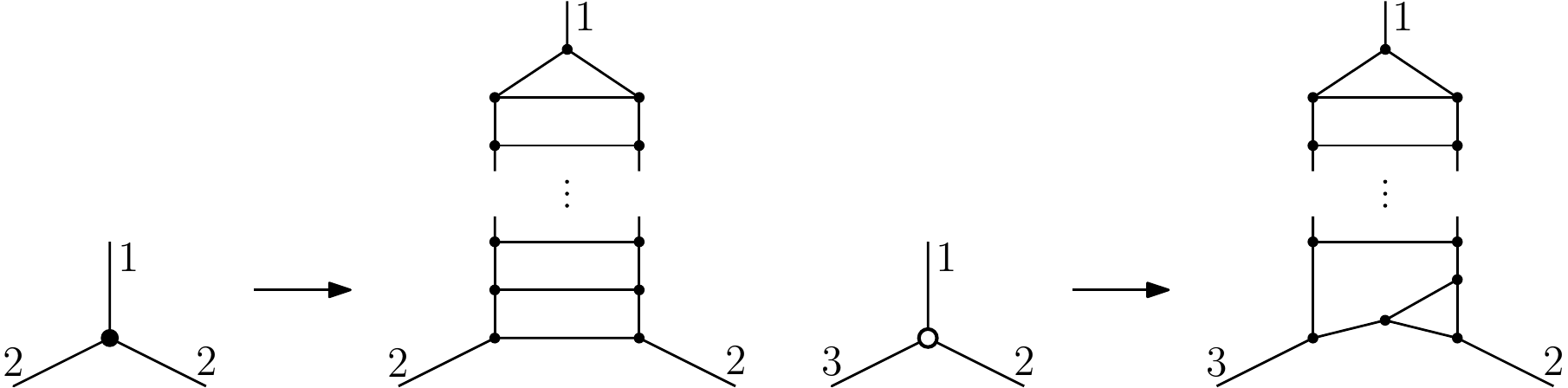}\\[1mm]
Figure 2: Operations $A$ (left-hand side) and $B$ (right-hand side).
\end{center}

Using operations $A$ and $B$, replace in $H$ each black vertex with a copy of $A_{r+2}$ and each white vertex with a copy of $B_r$, respecting the orientations given in Fig.~1 by the numbers 1, 2, 3. We obtain a planar graph $G$ that is clearly 3-connected and cubic. The circumference of $A_{r+2}$ and $B_r$ is $2r+5$. By construction, any cycle in $G$ of length greater than $2r+5$ has length at least $4r+15$, which is the length of the cycle bounding the face $F$ and also of the cycle bounding the face $F'$. Thus, for every $\ell \in \{ 2r+6, \dots, 4r+14 \}$, the graph $G$ contains no cycle of length $\ell$. Setting $k := 2r+6$, the proof is complete, since $G$ clearly has circumference at least $k$. \hfill $\Box$

\bigskip

Merker proves in~\cite{Me21} that for every integer $k \ge 4$ there exists a 3-connected cubic planar graph $G$ of circumference at least~$k$ which satisfies ${\cal C}(G) \cap [k,2k+1] = \emptyset$. In order to illustrate the construction yielding this result, Merker provides an example in~\cite[Fig.~2]{Me21}, which we will call $G$. We point out that, despite indeed explaining the construction method, this graph $G$ is not well chosen: $G$ does not satisfy the conditions Merker himself sets out and in consequence, there exists no positive integer $k \le |V(G)|$ ($G$ is hamiltonian) such that ${\cal C}(G) \cap [k, 2k+1] = \emptyset$. However, his proof is correct, only that $n$ (as defined in Merker's proof) must be chosen large enough in relation to $k$, as he himself states.

\bigskip

\bigskip

\noindent \textbf{Acknowledgements.} I thank Nico Van Cleemput for comments which improved the presentation of the above results. My research was supported by a Postdoctoral Fellowship of the Research Foundation Flanders (FWO).

\end{document}